\documentclass{amsart}

\usepackage{amsmath}
\usepackage{amsfonts}
\usepackage{amssymb}
\usepackage{graphicx}
\usepackage{cases}
\usepackage{hyperref}

\newtheorem{definition}{Definition}[section]
\newtheorem{theorem}[definition]{Theorem}
\newtheorem{lemma}[definition]{Lemma}
\newtheorem{example}[definition]{Example}
\newtheorem{corollary}[definition]{Corollary}
\newtheorem{notation}[definition]{Notation}
\theoremstyle{plain}

\newtheorem{notification}{Notification}[subsection]

\numberwithin{equation}{section}

\begin{document}
\pagestyle{plain}

\title{RELATIONSHIPS AND ALGORITHM IN ORDER TO ACHIEVE THE LARGEST PRIMES }

\author{\scriptsize Ali Zalnezhad}
\address{ Noshiravani University of Technology, Babol, Iran}
\email{ali.zalnejad@stu.nit.ac.ir }

\author{ Ghasem Shabani}
\address{University of Tabriz, Tabriz, Iran}
\email{ghasemshabani@ymail.com}

\author{Hossein Zalnezhad}
\address{Science and Research Branch of Islamic Azad University, Tehran, Iran}
\email{h.zalnejad@gmail.com}

\author{Mehdi Zalnezhad}
\address{ Noshiravani University of Technology, Babol, Iran}
\email{m.zalnezhad.313@gmail.com}

\keywords{Generalization the Mersenne's  theorem, Relations of Prime numbers, Algorithm}

\begin{abstract}
Today , prime numbers attained exceptional situation in the area of numbers theory and cryptography. As we know, the trend for accessing to the largest prime numbers due to using Mersenne theorem, although resulted in vast development of related numbers, however it has reduced the speed of accessing to prime numbers from one to five years. This paper could attain to theorems that are more extended than Mersenne theorem with accelerating the speed of accessing to prime numbers.

Since that time, the reason for frequently using Mersenne theorem was that no one could find an efficient formula for accessing to the largest prime numbers. This paper provided some relations for prime numbers that one could define several formulas for attaining prime numbers in any interval; therefore, according to flexibility of these relations, it could be found a new branch in the field of accessing to great prime numbers followed by providing an algorithm at the end of this paper for finding the largest prime numbers.

\end{abstract}
\maketitle

\section{Introduction}

Due to the importance of the primes, the Mathematicians have been investigating about them since long centuries ago. in 1801, Carl  Gauss, one of the greatest mathematician, submitted that the problem of distinguishing  the primes among the non-primes has been one of the outstanding problems of arithmetic \cite{1}.
Proving the infinity of prime numbers by Euclid is one of the first and most brilliant works of the human being in the numbers theory \cite{2}. Greek people knew prime numbers and were aware of their role as building blocks of other numbers. More, the most natural question asked by human being was this what order prime numbers are following and how one could find prime numbers? Until this time, there have been more attempts for finding a formula producing the prime numbers and or a model for appearance of prime numbers among other numbers and although they could be more helpful for developing the numbers theory, however, the complicated structure of prime numbers could not be decoded. During last years, the prime numbers attained an exceptional situation in the field of coding. For example, “RSA” system is one of the most applicable system in this field used in industries relying on prime numbers. “RSA” system is used in most computerized systems and counted as main protocol for secure internet connections used by states and huge companies and universities in most computerized systems \cite{3}. On 2004, Manindra Agrawal and his students in Indian Institute of Technology Kanpur could develop an algorithm called AKS for detecting prime numbers \cite{4}. On 2006, 2008, 2009 and recently on 2013, mathematics students in a project called detecting the Mersenne Prime Numbers by Computer Network GIMPS succeeded to discover the greatest prime number. All such cases indicate the importance of Mersenne theorem or any other approach for finding the largest prime numbers \cite{5}. Generalizing the Mersenne theorem, this paper could accelerate finding the largest prime numbers. in addition, there have been provided new equations and algorithm for attaining the largest primes.

\section{Generalization the Mersenne's  theorem }

\begin{definition}
Assume that $n$ is a natural number greater than 1, $Z_n$  related to n and natural numbers $a$ and $C$ are defined as below:
\begin{equation}
Z_n=\frac{(a+C)^n-a^n}{C}      \hspace{1cm}    (a \in N , n \in N , n > 1  , C \in N)  \nonumber \\  
\end{equation}
\end{definition}
\begin{theorem}
 If  $Z_n$  is a prime number, then $n$  is a prime number, too. 
 \end{theorem}
\begin{proof}
 If $n$ is not the prime number so we can write $n$ as the multiplication of two  natural numbers except $1$. meaning:
\begin{equation}
n=rs, r\not=1 , s\not=1          \nonumber \\
\end{equation}
\begin{equation}
Z_n=\frac{〖(a+c)〗^n-(a)^n}{c}=\frac{〖(a+c)〗^{rs}-(a)^{rs}}{c}=\frac{((a+c)^r )^s-(a^r )^s }{c}   \nonumber \\      
\end{equation}
\begin{equation}
=\frac{((a+c)^r-a^r )((a+c)^{n-r}+(a+c)^{n-2r} a^r+\ldots+a^{n-r})}{c}     \nonumber \\    
\end{equation}
\begin{equation}
= \frac{(a^r+{ra^{r-1} c}+\ldots+c^r-a^r )((a+c)^{n-r}+(a+c)^{n-2r} a^r+\ldots+a^{n-r})}{c}  \nonumber \\        
\end{equation}
\begin{equation}
=\frac{({ra^{r-1} c}+\ldots+c^r )((a+c)^{n-r}+(a+c)^{n-2r} a^r+\ldots+a^{n-r})}{c}       \nonumber \\
\end{equation}
\begin{equation}
=({ra^{r-1}}+\ldots+c^{r-1} )((a+c)^{n-r}+(a+c)^{n-2r} a^r+\ldots+a^{n-r})        \nonumber \\ 
\end{equation}
Therefore, $Z_n$ is not the prime number. so, $n$ must be a prime number.
\end{proof}
\subsection*{Note}
 This theorem is a generalization for Mersenne theorem in which  $a$ and $C$ are arbitrary natural numbers.
\begin{lemma}
  If in the theorem $2.2$, C is chosen as a multiple to $a$  and  $a \geq 2$,  thus, $Z_n$ will not be a prime number.
  \end{lemma}
\begin{proof}
Suppose:
\begin{equation}
C=ad  \hspace{0.5cm} \{ a \in N , a \geq 2 ,  d \in N \}       \nonumber \\
\end{equation}
Therefore:
\begin{equation}
Z_n=\frac{(a+C)^n-a^n}{C} = \frac{(a+ad)^n-a^n}{ad} = \frac{a^n (1+d)^n-a^n}{ad}      \nonumber \\
\end{equation}
\begin{equation}
=\frac{〖a^{n-1} (1+d)〗^n-a^{n-1}}{d}  =\frac {a^{n-1} ((1+d)^n-1)}{d}      \nonumber \\
\end{equation}
\begin{equation}
=\frac{a^{n-1}((d+1)^1-1^1 ))((d+1)^{n-1}+(d+1)^{n-2}+\ldots+1)}{d}      \nonumber \\
\end{equation}
\begin{equation}
=\frac{a^{n-1}(d)((d+1)^{n-1}+(d+1)^{n-2}+\ldots+1)}{d}     \nonumber \\
\end{equation}
\begin{equation}
=a^{n-1} ((d+1)^{n-1}+(d+1)^{n-2}+\ldots+1)     \nonumber \\
\end{equation}
\begin{equation}
=a^{n-1} k   \hspace{1cm}  (a \not=1 ,k \not=1)   \nonumber \\
\end{equation}
The last equality shows that  $Z_n$ is not a prime number.
\end{proof}
\section{A Specific State of Generalized Theorem }
\begin{definition}
Suppose $n$ is a natural number greater than $1$, function $Z_n$ related to $n$ and natural number $a$ are defined as below:
\begin{equation}
Z_n=〖(a+1)〗^n-a^n    \hspace{1cm}  (a \in N , n \in N ,n>1)  \nonumber \\
\end{equation}
\end{definition}
\begin{theorem}
If $Z_n$ is a prime number, then $n$ is a prime number, too.
\end{theorem}
\begin{lemma}
In this theorem  $3.2 $, based on $n$ constant, please  consider a sequence
\end{lemma}
\begin{equation}
  Z_a=〖(a+1)〗^n-a^n  \nonumber \\
\end{equation}
We prove that sequence ${\{Z_a\}}_{a \in N}$ is strictly ascending, i.e.
\begin{equation}
 Z_a>Z_{a-1} \hspace{1cm}  (a \in N)  \nonumber \\
\end{equation}
To prove the last inequality, we write:
\begin{equation}
(〖a+1)〗^n-(a)^n>(a)^n-(a-1)^n  \nonumber \\
\end{equation}
\begin{equation}
〖(a+1)〗^n+〖(a-1)〗^n>2〖(a)〗^n  \nonumber \\
\end{equation}
\begin{multline}
a^n+\frac{na^{n-1}}{1!}+\frac{n(n-1) a^{n-2}}{2!}+\ldots+na+1 \\ +a^n-\frac{na^{n-1}}{1!}+\frac{n(n-1) a^{n-2}}{2!}+\ldots+(-1)^{n-1} na+〖(-1)〗^n >2(a)^n  \nonumber \\
\end{multline}
\begin{equation}
2(\frac{n(n-1) a^{n-2}}{2!}+\ldots+na\frac{1+(-1)^{n-1} }{2}+\frac{1+(-1)^n }{2})>0  \nonumber \\
\end{equation}
Status 1. If $n$ is a multiple of $2$:
\begin{equation}
2(\frac{n(n-1) a^{n-2}}{2!}+\ldots+0+1)>0  \nonumber \\
\end{equation}
Status 2. If $n$ is not a multiple of $2$:
\begin{equation}
2(\frac{n(n-1) a^{n-2}}{2!}+\ldots+na+0)>0  \nonumber \\
\end{equation}
Therefore, inequity is accepted.
\begin{corollary}
In this theorem, each number is higher than Mersenne number, meaning:
\begin{equation}
 Z_a>Z_1 \rightarrow (a+1)^n-(a)^n>2^n-1 \hspace{1cm}  (a \in N,a>1 ,n \in N ,n>1)   \nonumber \\
\end{equation}
\end{corollary}
\section{Obtaining ALL Primes In Given Interval}
\begin{definition}
Suppose $a$  be a natural number  and  $2,C_1,C_2,\ldots,C_n$ are the primes smaller than or equal $\sqrt{a}$   and $k$ , $m$ are natural numbers which limitations are intended for them indicated As follows:  
\begin{equation}
  K\not= 〖C 〗_i m+  \frac{〖C 〗_i-1}{2} \hspace{0.3cm} for \hspace{0.3cm} \frac{〖C 〗_i-1}{2} \leq m<\frac{a-C_i}{2C_i } \hspace{0.5cm}  (i=1,2,3,\ldots,n)  ,  m \in N    \nonumber \\
\end{equation}
\end{definition}
\begin{theorem}
Assume that $R$  is a function of $ K$  which is displayed as bellow:
\end{theorem}
\begin{equation}
 R=2K+1  \hspace{0.5cm}  \{ K \in N ,  K< \frac{a-1}{2} \}  \nonumber \\
\end{equation}
If the $K$ and $m$ circumstances are followed, $R$ can obtain all the primes less than $a$. 
\begin{proof}
Knowing that $R$ is odd, because it is non  prime, therefore it comprises from two odd numbers except $1$, and because $R< a$, $R$ has at least a prime factor $ \leq  \sqrt{a}$.
Therefore, $R$ is divided at least on one of the prime factors ${\{C_ i\}}_{i=1,2,\ldots,n}$ .
\begin{equation}
 R=2K+1  \hspace{0.25cm}, \hspace{0.25cm}  R=C_i (2m+1)  \nonumber \\
\end{equation}
\begin{equation}
 C_i (2m+1)=2K+1 \rightarrow K=〖C 〗_i m+  \frac {〖C 〗_i-1}{2}  \hspace{0.25cm} (i=1,2,\ldots,n) \nonumber \\
\end{equation}
It is clear that above equalities are in discrepancy of the assumption of the theorem.
\end{proof}
\subsection{The reason of putting intervals m}
\begin{enumerate}

 \item if :
\begin{equation}
\exists i \in N :m \geq \frac{a-C_i}{2C_i } \rightarrow C_i (2m+1) \geq a \rightarrow R\geq a. \nonumber \\
\end{equation}
\item interval $\frac{〖C 〗_ i-1}{2} \leq m$:

It is clear that by putting minimum $m$ in the definition $4.1$ minimum $ k$  followed by minimum  $R$  is obtained as below:
\begin{equation}
m=\frac{C_i-1}{2} \rightarrow K=  \frac{〖C 〗_i-1}{2 }(C_i+1) \rightarrow R=〖C_i〗^2 \nonumber \\
\end{equation}
\end{enumerate}
According to recent equation, it is obvious that being as prime number in prime numbers smaller  than  $〖C_i〗^2$,   R may not be divided into prime factors smaller than $C_i$. On the other hand, it is not necessary to see if prime numbers smaller than  $〖C_i〗^2$   are divided into $〖C 〗_i $ to detect it as a prime number. Indeed, for obtaining the prime numbers, we only require $ R$ in ${C_ i}^2 \leq R<a $ to enter the provision of prime factor    $C_ i  $ .
\subsection*{Note}
If $C_{n+1}$ is considered as a prime number bigger than $C_n$,   we could use $〖C_{n+1 }〗^2-1$ instead of $a$ in this theorem because prime numbers smaller than $\sqrt{a}$ include prime numbers smaller than $\sqrt{{C_{n+1}}^2-1}$ .
\begin{example}
Prime numbers smaller than 120: 
\begin{equation}
R=2K+1  \hspace{0.5cm}  K \in N \hspace{0.25cm} , \hspace{0.25cm}   K< \frac{a-1}{2} \hspace{0.25cm}   or \hspace{0.25cm} K<\frac{〖C_{n+1}〗^2-1}{2} \nonumber \\
\end{equation}
\begin{center}
\{ $ a \in N$ \hspace{0.25cm},\hspace{0.25cm}   Prime numbers smaller than $\sqrt{a}$:   ${\{2, C _1 ,C _2  ,C_3\} }$ \} \nonumber \\
\end{center}
\begin{equation}
K\not= 〖C 〗_i m+  \frac{ 〖C 〗_i-1}{2} \hspace{0.25cm} for \hspace{0.25cm} \frac{〖C 〗_ i-1}{2} \leq m<\frac{a-C_i}{2C_i } \hspace{0.25cm}  (i=1,2,3) \nonumber \\
\end{equation}
\begin{equation}
a=120 , C_1=3 , C_2=5 , C_3=7  \nonumber \\
\end{equation}
\begin{equation}
R=2K+1   \hspace{0.25cm}  K \in N ,  K<\frac{120-1}{2}    \rightarrow K<59.5  \nonumber \\
\end{equation}
\begin{equation}
K\not= 3m+  \frac{ 3-1}{2} \hspace{0.2cm} for \hspace{0.2cm} \frac{ 3-1}{2}\leq m<\frac{120-3}{2\times3}   \rightarrow   K\not= 3m+ 1  \hspace{0.25cm} for  \hspace{0.25cm} 1\leq m<\frac{117}{6}  \nonumber \\
\end{equation}
\begin{equation}
K\not= 4,7,10,13,16,19,22,25,28,31,34,37,40,43,46,49,52,55,58  \nonumber \\
\end{equation}
\begin{equation}
K\not= 5m+  \frac{ 5-1}{2} \hspace{0.2cm} for \hspace{0.2cm} \frac{ 5-1}{2}\leq m<\frac{120-5}{2\times5}   \rightarrow   K\not= 5m+ 2  \hspace{0.25cm} for  \hspace{0.25cm} 2\leq m<\frac{115}{10}  \nonumber \\
\end{equation}
\begin{equation}
K\not= 12,17,22,27,32,37,42,47,52,57  \nonumber \\
\end{equation}
\begin{equation}
K\not= 7m+  \frac{ 7-1}{2} \hspace{0.2cm} for \hspace{0.2cm} \frac{ 7-1}{2}\leq m<\frac{120-7}{2\times7}   \rightarrow   K\not= 7m+ 3  \hspace{0.25cm} for  \hspace{0.25cm} 3\leq m<\frac{113}{14}  \nonumber \\
\end{equation}
\begin{equation}
K\not= 24, 31, 38, 45, 52, 59  \nonumber \\
\end{equation}

\begin{equation}
R=3,5,7,11,13,17,19,23,29,\ldots,113  \nonumber \\
\end{equation}
\end{example}
\section{Relation 1: Determination of Primes less than a given natural number}
\begin{definition}
Suppose  $a$ be the natural number and $C_1,C_2,\ldots,C_n$ are the primes smaller than or equal $\sqrt{a}$  and also consider that $C_{n+1},C_{n+2},\ldots$ are the primes larger than $\sqrt{a}$ . Suppose that $b_1,b_2$  and $K$ be the members of the natural numbers and also   $m_1,m_2,\ldots$ be the members of the account numbers, these variables are selected arbitrarily. Function $R$  related to values ${\{C_i\}}_{i \in N}$, natural numbers $K, b_2, b_1$ and arithmetic numbers ${\{m_i\}}_{i \in N}$ are defined as below:
\begin{equation}
R=(-1)^{b_1 } \times K \times C _1  \times C_2   \times \ldots \times  C_n + (-1)^{b_2 } \times {C_{n+1}} ^{m_1} \times {C_{n+2}}^{m_2} \times \ldots 
  \nonumber \\
\end{equation}
\begin{equation}
(K \in N \hspace{0.25cm}  b_1,b_2 \in N    \hspace{0.25cm} m_1  ,m_2 , \ldots \in W) \nonumber \\
\end{equation}
\end{definition}
\begin{theorem}
If $ R$ be as the natural number less than $a$, then $R$ is a prime number.
\end{theorem}
 \begin {proof}
 If $R$  is not prime number, it has a prime factor $ \leq \sqrt{R}$. On the other side, because $ R<a$, $R$ has at least one prime factor $ \leq \sqrt{a}$. So, it is arbitrarily supposed that $R$ is divisible in $\{C_i\}_{i =1,2.\ldots,n}$.
\begin{equation}
\exists i: \frac{R}{C_i} \in N  \nonumber \\
\end{equation}
\begin{equation}
\frac{R}{C_i}=\frac{(-1)^{b_1 } \times K \times C _1  \times C_2 \times \ldots \times  C_n}{C_i} +\frac{ (-1)^{b_2 } \times {C_{n+1}} ^{m_1} \times {C_{n+2}}^{m_2} \times \ldots}{C_i}   \nonumber \\
\end{equation}

Because $C_i$ is not denominator of any $\{C_j\}_{j = n+1,n+2,\ldots}$ . we have:
\begin{equation}
\frac{R}{C_i}=\pm natural \hspace{0.1cm}  number \pm decimal  \hspace{0.1cm}  number \nonumber \\
\end{equation}
\begin{equation}
\frac{R}{C_i}= decimal  \hspace{0.1cm}  number \nonumber \\
\end{equation}
We reached a contradiction to the assumption. thus, the theorem was verified.  
\end{proof}
\begin{example}
We obtain prime numbers smaller than 119.\\
\begin{equation*}
\begin{cases}
\sqrt{119}\approx10 \\
\text{Primes smaller than or equal} \hspace{0.1cm} 10=2 ,3 ,5 ,7  \\
 
\text{ Primes larger than} \hspace{0.1cm} 10=11 ,13 ,17, \ldots  
\end{cases}\\
\end{equation*}
\begin{equation}
R=(-1)^{b_1 } \times K \times 2  \times 3 \times 5 \times 7  + (-1)^{b_2 } \times 11 ^{m_1} \times 13^{m_2} \times 17^{m_3 } \times \ldots 
  \nonumber \\
\end{equation}
\begin{equation}
R=(-1)^{b_1 } \times K \times 210  + (-1)^{b_2 } \times 11 ^{m_1} \times 13^{m_2} \times 17^{m_3 } \times \ldots    \nonumber \\
\end{equation}
\end{example}
\begin{center}
Table 1

\begin{tabular}{|c|c|c|c|c|c|}\hline
$b_1$ & 2 &1 & 2 & 1& 1  \\ \hline
$b_2$ &1 &2 &1 &2 &2  \\ \hline
$K$ &1 &6 &1 &9 &7  \\ \hline
$m_i$ & $m_1=2$ & $m_1=2$  & $m_1=m_2=1 $ & $m_1=1 ,m_2=2 $  & $m_1=2 , m_2=1 $  \\
 & $m_{n}=0 $ &  $m_{n}=0 $ &  $ m_{n}=0 $ &  $ m_{n}=0$ &  $ m_{n}=0 $ \\  
&  $ n \geq 2 $ & $ n \geq 2 $ & $ n \geq 3 $ &  $ n \geq 3 $ & $ n \geq 3  $ \\  \hline
$R$ &89 &71 &67 &31 &103  \\ \hline
\end{tabular}
   
 And continuing so \ldots

\end{center}

\subsection{The sensible hints about the relation rang}
Suppose that: $ C_1<C_2<\ldots<C_n<C_{n+1}<\ldots $ (the order is considered in the primes)
\begin{notification}
General speaking, the theorem $ 5.2 $ comes true to   $ 1<R \leq {C_{n+1}}^2-1  $   because $ R $  includes the same primes.
\end{notification}
\begin{notification}
$R$ is obviously not divisible to $C_1,C_2,\ldots ,C_n$  and According to prime of the number $R$, we have:
\begin{equation}
C_n<R\leq {C_{n+1}}^2-1 \nonumber \\
\end{equation}
\end{notification}
\begin{notification}
 If the arithmetic sequence $\{{m_i}\}_{i\in N}$  is descending, i.e:
$ m_1 \geq m_2 \geq m_3 \geq \ldots $
Then: $C_n<R \leq {C_{the \hspace{0.1cm} first \hspace{0.1cm} m \hspace{0.1cm} which \hspace{0.1cm} turned \hspace{0.1cm} into \hspace{0.1cm} zero}}^2-1$
\end{notification}
\subsubsection*{Example}
If $m_4=0 ,{\{m_i \not=0\}}_{ i=1,2,3}$ then $C_n <R \leq {C_{n+4 }}^2-1 $
\begin{notification}
To attain prime numbers, we divide the intervals as below:
\begin{equation}
\{ {C_{n+i}}^2,{C_{n+i+1}}^2-1\} \hspace{0.2cm} , \hspace{0.2cm}  i\in N \nonumber \\
\end{equation}
With regard to the relationship easier to be written. In example of the primes less than $100$, the rang can be divided into three sections of $(4-24)$, $(25-48)$ and $(49-100)$. Then, a distinct relation asserted for each.
\end{notification} 
\subsubsection*{Example}
Prime numbers smaller than 48 :
\begin{equation}
(1-24)   \rightarrow R=6K \pm 1 \nonumber \\
\end{equation}
\begin{equation}
K=1  \rightarrow R=(5,7 ) \hspace{0.35cm}  K=3  \rightarrow R=(17,19) \hspace{0.35cm}  K=4  \rightarrow R=23    \nonumber \\
\end{equation}
\begin{equation}
(25-48)   \rightarrow R=(-1)^{b_1 } \times K \times 30  + (-1)^{b_2 } \times 11 ^{m_1} \times 13^{m_2} \times 17^{m_3 } \times \ldots    \nonumber \\
\end{equation}
\begin{equation}
R=30 \pm 1=(29,31) , R=30 \pm 7=(23,37) , R=30 \pm 13=(17,43)    \nonumber \\
\end{equation}
\begin{center}
 And continuing so \ldots
\end{center}
\subsection{Special case of relation 1}
If [$\sqrt{a}$ ] is considered to be equal d $\rightarrow$ ([$\sqrt{a}$ ]=d), and then substituted in k, i.e: filling blank space between the primes in the  relation 1, we can conclude:
\begin{equation}
 K \times C _1  \times C_2 \times C _3  \times \ldots \times  C_n =k_1 \times d! \hspace{0.5cm} (k_1 \in N)  \nonumber \\
\end{equation}
\hspace{7cm}$ \rightarrow$
\begin{equation}
 R=(-1)^{b_1 } \times k_1 \times d!  + (-1)^{b_2 } \times {C_{n+1}} ^{m_1} \times {C_{n+2}}^{m_2} \times {C_{n+3}}^{m_3 } \times \ldots \nonumber \\
\end{equation}
\begin{equation}
\{ (b_1,b_2) \in N   \hspace{0.5cm}  k_1 \in N   \hspace{0.5cm}  m_i \in W \} \nonumber \\
\end{equation}
Also, if  $  (m_2,m_3,\ldots ,m_n=0)  \rightarrow R=(-1)^{b_1}×k_1×d!+(-1)^{b_2 }×{C_{n+1}}^{ m_1}$ .
Therefore, value of $ R$ will be a prime number by the provision of $ C_n< R \leq {C_{n+1}}^2-1.$
(of course this relation is utilized for formula simplicity)
\begin{notation}
 We could also replace the known prime numbers in the relation 1 and replace the unknown ones as a multiple of K; by this, we can attain a simpler equation.
 \end{notation}
\subsection*{Note}
If prime number is very big, one can take assistance from Notification 5.1.4 and Notation 5.4 .
\section{Relation 2: Determination of Primes less than a given natural number}
\begin{definition}
Suppose  $a$ be the natural number  and  $C_z  ,C_x,\ldots,C_v  ,C_m,C_j, \ldots ,C_n$ are the primes smaller than or equal $\sqrt{a}$ and also consider that $C_{n+1}$ is prime number larger than $\sqrt{a}$ . $R$ functions are defined as follows:
\begin{equation}
R=(-1)^{b_1} \times C _z  \times  C _x \times \ldots \times C _v \times k_1+(-1)^{b_2 } \times C_m \times C_ j \times \ldots \times C_n \times k_2 \nonumber \\
\end{equation}
\begin{equation}
+ (-1)^{b_3} \times C_z  \times C _x \times \ldots \times C _v \times C_m \times C_j \times \ldots \times C_n \times k_3  \nonumber \\
\end{equation}
\begin{equation}
(b_1 , b_2 , b_3 ,  k_1 , k_2 \in N,k_3 \in W)  \nonumber \\
\end{equation}
In addition, $k_1$ is not divisible on any of prime numbers $ C_m,C_j,\ldots ,C_n $ and $k_2$ also is not divisible on any prime numbers $ C_z  ,C _x  , \ldots ,C _v $ .
\end{definition}
\subsubsection{Notification}
Function $R$ has been comprised from three sections. Section one includes a part of the primes less or equal $\sqrt{a}$ ,  section two comprises the primes which are not contained in section one, and the third section includes all the primes less or equal $\sqrt{a}$ .
\begin{theorem}
If  $R$ be the natural number less than $a$, then $R$ is a prime number.
\end{theorem}
\begin{proof}
If  $R$ is not a Prime number ,it has a prime factor  $ \leq \sqrt{R}$. On the other side, because $R<a, R$ has at least one prime factor $\leq \sqrt{a}$. so, it is arbitrarily supposed that $R$ is divisible in $\{C_{\beta}\}_{\beta=z,x,\ldots,v,m,j,\ldots,n}$.
\begin{equation}
\exists \beta:  \frac{R}{C_{\beta}} \in N  \nonumber \\
\end{equation}
\begin{equation}
\frac{R}{C_{\beta}}=\frac{(-1)^{b_1} \times C _z  \times  C _x \times \ldots \times C _v \times k_1}{C_{\beta}}+\frac{(-1)^{b_2 } \times C_m \times C_ j \times \ldots \times C_n \times k_2}{C_{\beta}}  \nonumber \\
\end{equation}
\begin{equation}
+ \frac{(-1)^{b_3} \times C_z  \times C _x \times \ldots \times C _v \times C_m \times C_j \times \ldots \times C_n \times k_3}{C_{\beta}}  \nonumber \\
\end{equation}
$C_{\beta}$ may not be located in one of two parts first or second. therefore:
\begin{equation}
\frac{R}{C_{\beta}}=\pm natural \hspace{0.1cm} number \pm decimal \hspace{0.1cm} number \pm natural \hspace{0.1cm} number   \nonumber \\
\end{equation}
\begin{equation}
 \frac{R}{C_{\beta}} =decimal \hspace{0.1cm} number   \nonumber \\
\end{equation}
That is inconsistent with being \begin{large} $\frac{R}{C_{\beta}}$  \end{large} as natural number. Therefore, theorem has been proved.
\end{proof}
\begin{example}
Primes less than 119: 
\begin{equation*}
\begin{cases}
\sqrt{119}\approx10 \\
\text{Primes smaller than or equal}  \hspace{0.1cm} 10=2 ,3 ,5 ,7  \\
\end{cases}\\
\end{equation*}
In order to simplify the relation, the amount of $k_3$ can be ascertained as zero. Then, we can have the following relation $ \rightarrow k_3=0$.
\begin{equation}
R=(-1)^{b_1} \times 2  \times  7 \times k_1+(-1)^{b_2 } \times 3 \times 5 \times k_2 \nonumber \\
\end{equation}
\begin{equation}
R=(-1)^{b_1} \times 14 \times k_1+(-1)^{b_2 } \times 15 \times k_2 \nonumber \\
\end{equation}
\begin{equation}
(b_1 , b_2 \in N   \hspace{0.5cm}  \frac{ K_1}{3}\notin N , \frac{K_1}{5} \notin N \hspace{0.5cm}   \frac{K_2}{2} \notin N , \frac{K_2}{7} \notin N )                         \nonumber \\
\end{equation}

\begin{center}
Table 2

\begin{tabular}{|c|c|c|c|c|c|c|c|}\hline

$b_1$ & 2 &2 & 2 & 2& 2&1&1  \\ \hline
$b_2$ & 2 &2 & 2 & 1& 2&2&2  \\ \hline
$K_1$ & 1 &2 & 1 & 11 & 2&2&1  \\ \hline
$K_2$ & 3 &3 & 5 & 5& 1&3&5  \\ \hline
$R$ & 59 &73 & 89 & 79 & 43&17&61  \\ \hline
\end{tabular}

And continuing so \ldots
\end{center}
\end{example}
\subsection{The sensible hints about the relation rang}
Assume that $C_n$ is the biggest prime number $ \leq \sqrt{a}$ .
\begin{notification}
General speaking,the theorem $ 6.2 $ comes true to   $ 1<R \leq {C_{n+1}}^2-1  $   because $ R $  includes the same primes.
\end{notification}
\begin{notification}
$R$ is obviously not divisible to $C_z  ,C _x ,\ldots,C_v  ,C_m,C_j,\ldots,C_n$  and According to prime of the number $R$, we have:
\begin{equation}
C_n<R\leq {C_{n+1}}^2-1 \nonumber \\
\end{equation}
\end{notification}
\begin{notification}
To attain prime numbers, we divide the intervals as below:
\begin{equation}
\{ {C_{n+i}}^2,{C_{n+i+1}}^2-1\} \hspace{0.2cm} , \hspace{0.2cm}  i\in N \nonumber \\
\end{equation}
With regard to the relationship easier to be written.
In example of the primes less than $100$, the rang can be divided into three sections of $(4-24)$, $(25-48)$ and $(49-100)$. Then, a distinct relation asserted for each. 
\end{notification}
\subsubsection*{Example}
Prime numbers smaller than 48 :
\begin{equation}
(1-24)   \rightarrow R=(-1)^{b_1 }\times 3 \times k_1+(-1)^{b_2} \times 2 \times k_2    \hspace{0.5cm} ( \frac{ K_1}{2} \notin N  ,  \frac{ K_2}{3} \notin N ) \nonumber \\
\end{equation}
\begin{equation}
(25-48)   \rightarrow R=(-1)^{b_1 }\times 3 \times 2 \times k_1+(-1)^{b_2} \times 5 \times k_2    \hspace{0.3cm} ( \frac{ K_1}{5} \notin N  ,  \frac{ K_2}{2} \notin N  ,  \frac{ K_2}{3} \notin N ) \nonumber \\
\end{equation}
\section{Relation 3: Determination of Primes less than a given natural number}
\begin{definition}
Suppose $a$ be the natural number and $C_1  , C _2 , C_3 ,\ldots, C_{n-2} , C_{n-1} , C_n $ are primes smaller than or equal $\sqrt{a}$ and also consider that $C_{n+1}$ be the prime number larger than $\sqrt{a}$ and $k_1 , k_2 , \ldots , k_n , k_{n+1} , b_1 , b_2 , \ldots , b_n , b_{n+1}$    be the members of the account numbers. Now, we display $R$ function as below:
\begin{equation}
R=(-1)^{b_1} \times C_n  \times C _{n-1} \times \ldots \times C_3 \times C_2 \times k_1+(-1)^{b_2} \times C_n \times C_{n-1} \times  \ldots \times C_3 \times C_1 \times k_2   \nonumber \\
\end{equation}
\begin{equation}
+\ldots+(-1)^{b_n} \times C_{n-1}  \times C_{n-2}\times \ldots \times C_3 \times C_2 \times C_1 \times k_n  \nonumber \\
\end{equation}
\begin{equation}
+(-1)^{b_{n+1} } \times C_n \times C_{n-1}  \times C_{n-2} \times \ldots \times C_3 \times C_2 \times C_1 \times k_{n+1}  \nonumber \\
\end{equation}
\begin{equation}
\{ b_1,b_2,…,b_n,b_{n+1} \in N   \hspace{0.5cm}   k_1,k_2 , \ldots , k_n \in N , k_{n+1} \in W \}  \nonumber \\
\end{equation}
\begin{equation}
 \{ \frac{ k_1}{C_1} \notin N , \frac{ k_2}{C_2} \notin N , \ldots , \frac{ k_n}{C_n} \notin N \} \nonumber \\
\end{equation}
\end{definition}
\begin{theorem}
If  $R$ be the natural number less than $a$, then $R$ is a prime number.
\end{theorem}
\begin{proof}
If $R$ is not a Prime number ,it has a prime factor $ \leq \sqrt{R}$. On the other side, since $ R<a, R $ has at least one prime factor $ \leq \sqrt{a}$.
So, it is arbitrarily supposed that $R$ is divisible in $\{C_i\}_{i =1,2.\ldots,n}$ .
\begin{equation}
\exists i: \frac{R}{C_i} \in N  \nonumber \\
\end{equation}
\begin{equation}
\frac {R}{C_i} = \frac{(-1)^{b_1} \times C_n  \times C _{n-1} \times \ldots \times C_3 \times C_2 \times k_1}{C_i} \nonumber \\
\end{equation}
\begin{equation}
+ \frac{(-1)^{b_2} \times C_n \times C_{n-1} \times  \ldots \times C_3 \times C_1 \times k_2}{C_i}   \nonumber \\
\end{equation}
\begin{equation}
+\ldots+ \frac{(-1)^{b_n} \times C_{n-1}  \times C_{n-2}\times \ldots \times C_3 \times C_2 \times C_1 \times k_n}{C_i}  \nonumber \\
\end{equation}
\begin{equation}
+ \frac{(-1)^{b_{n+1} } \times C_n \times C_{n-1}  \times C_{n-2} \times \ldots \times C_3 \times C_2 \times C_1 \times k_{n+1}}{C_i}  \nonumber \\
\end{equation}
Because there is only one quotient that is not an integer, we have:
\begin{equation}
\frac{R}{C_i}=decimal \hspace{0.1cm} number    \nonumber \\
\end{equation} 
That is inconsistent with being \begin{large} $\frac{R}{C_i} $  \end{large} as natural number. Therefore, theorem has been proved.
\end{proof}

\begin{example}
Primes less than 119: 
\begin{equation*}
\begin{cases}
\sqrt{119}\approx10 \\
\text{Primes smaller than or equal}  \hspace{0.1cm} 10=2 ,3 ,5 ,7  \\
\end{cases}\\
\end{equation*}

In order to simplify the relation, the amount of $k_{n+1}$ can be ascertained as zero. Then, we can have the following relation $\rightarrow k_{n+1}=0$.

\begin{equation}
R=(-1)^{b_1} \times 7 \times 5 \times 3 \times k_1+(-1)^{b_2} \times 7 \times 5 \times 2 \times k_2 \nonumber \\
\end{equation}
\begin{equation}
+(-1)^{b_3} \times 7 \times 3 \times 2 \times k_3+(-1)^{b_4} \times 5 \times 3 \times 2 \times k_4  \nonumber \\
\end{equation}
\begin{equation}
R=(-1)^{b_1} \times 105 \times k_1+(-1)^{b_2} \times 70 \times k_2 \nonumber \\
\end{equation}
\begin{equation}
+(-1)^{b_3} \times 42 \times k_3+(-1)^{b_4} \times 30 \times k_4  \nonumber \\
\end{equation}
\begin{equation}
\{ b_1,b_2,b_3,b_4  \in N   \hspace{0.5cm}   k_1,k_2 ,k_3 , k_4 \in N \}  \nonumber \\
\end{equation}
\begin{equation}
 \{ \frac{ k_1}{2} \notin N , \frac{ k_2}{3} \notin N  , \frac{ k_3}{5} \notin N ,  \frac{ k_4}{7} \notin N \} \nonumber \\
\end{equation}

\begin{center}
Table 3

\begin{tabular}{|c|c|c|c|c|}\hline

$b_1$ & 2 &2 & 2 & 2 \\ \hline
$b_2$ & 1 &2 & 2 & 2  \\ \hline
$b_3$ & 2 &1 & 1 & 1  \\ \hline
$b_4$ & 2 &1 & 1 & 1  \\ \hline
$K_1$ & 1 &1 & 1 & 1   \\ \hline
$K_2$ & 1 &1 & 1 & 2 \\ \hline
$K_3$ & 1 &2 & 1 & 2 \\ \hline
$K_4$ & 1 &1 & 1 & 2 \\ \hline
$R$ & 107 &61 & 103 & 101  \\ \hline
\end{tabular}

And continuing so \ldots
\end{center}
\end{example}
\subsection{The sensible hints about the relation rang}
Suppose that: $ C_1<C_2<\ldots<C_n<C_{n+1}<\ldots $ (the order is considered in the primes)
\begin{notification}
General speaking, the theorem $ 7.2 $ comes true to   $ 1<R \leq C_{n+1}^2-1  $   because $ R $  includes the same primes.
\end{notification}
\begin{notification}
$R$ is obviously not divisible to $C_1,C_2,\ldots ,C_n$  and According to prime of the number $R$, we have:
\begin{equation}
C_n<R\leq {C_{n+1}}^2-1 \nonumber \\
\end{equation}
\end{notification}
\begin{notification}
To attain prime numbers, we divide the intervals as below:
\begin{equation}
\{{C_{n+i}}^2,{C_{n+i+1}}^2-1\} \hspace{0.2cm} , \hspace{0.2cm}  i\in N \nonumber \\
\end{equation}
With regard to the relationship easier to be written.
In example of the primes less than $100$, the rang can be divided into three sections of $(4-24)$, $(25-48)$ and $(49-100)$. Then, a distinct relation asserted for each. 
\end{notification}
\subsubsection*{Example}
Prime numbers smaller than 48 :
\begin{equation}
(1-24)   \rightarrow R=(-1)^{b_1 } \times 3 \times k_1+(-1)^{b_2} \times 2 \times k_2  \} \nonumber \\
\end{equation} 
\begin{equation}
\{ (b_1,b_2 \in N) ,  ( k_1, k_2  \in N) , (\frac{k_1}{2} \notin N , \frac{k_2}{3} \notin N \}   \nonumber \\
\end{equation}
\begin{equation}
(25-48)   \rightarrow R=(-1)^{b_1} \times 3 \times 5 \times k_1+(-1)^{b_2} \times 5 \times 2 \times k_2+(-1)^{b_3} \times 3 \times 2 \times k_3    \nonumber \\
\end{equation}
\begin{equation}
\{ (b_1,b_2,b_3 \in N) , ( k_1 , k_2 , k_3 \in N) , (\frac{k_1}{2} \notin N , \frac{k_2}{3} \notin N ,  \frac{k_3}{5} \notin N) \}   \nonumber \\
\end{equation}
\begin{center}
 And continuing so \ldots
\end{center}
\section{An Algorithm to obtain the Largest Prime Numbers}
By integrating the relations, particularly using the Relation 2 and Notation 5.4, we can attain an algorithm to obtain the largest prime number.; one of them is as below:
\subsection*{step 1}
Insert the largest known prime number in $a$: $a=2^{57,885,161}-1$
\subsection*{step 2}
Insert $(a-2)$ in $b$, meaning:  $b= a-2$
\subsection*{step 3}
Display the multiple of natural numbers $\leq b $ alternatively and as below and put it in $c$, i.e: $ c=b \times (b-2)  \times (b-4)  \times \ldots \times 9  \times 7  \times 5  \times 3  \times 1$
\subsection*{step 4}
Omit the Composite numbers in this product as far as you can and then calculate the multiplication operation.
\subsection*{Note}
 It is obvious that c includes all the primes less than or equal $\sqrt{a}$ except 2.
\subsection*{step 5}
First of all, assign $k$ as one  $(k=1)$ so as to obtain the minimum of $n$. then in the following equation, give a counter to $n$ through the minimum $n$ so long as k be a member of the natural numbers $(k \in N)$. meanwhile, $k$ should not be even.
\begin{equation}
\frac{a+2^n}{c} < k \leq \frac{a^2-1+2^n}{c}  \hspace{0.25cm} \{ n \in N  ,  \frac{k}{2} \notin N \}  \nonumber \\
\end{equation}
\subsection*{step 6}
Then insert the obtained $k$ and $n$ in the following equation to obtain the prime number . 
\begin{center}
$R=ck-2^n$
\end{center}
\subsection*{Note}
Steps 5 and 6 indicate that if $a<ck-2^n \leq a^2-1 \hspace{0.2cm}  ( \frac{ k}{2} \notin N ) , $   then $ck-2^n$ is a prime number.
\begin{proof}
If  $R$ is not prime number, it has a prime factor $ \leq \sqrt{R}$. On the other side, because $R<a $ , $R$ has at least one prime factor $ \leq \sqrt{a}.$  Therefore, with no interruption in the generality of the subject, we can assume that $R$ is divisible on $7$:                              
\begin{equation}
\frac{ck-2^n}{7} \in N   \nonumber \\
\end{equation}
\begin{equation}
\frac{ck-2^n}{7} =\frac{ b \times (b-2)  \times (b-4)  \times \ldots \times 9  \times 7  \times 5  \times 3  \times 1\times k}{7}-\frac{2^n}{7} \nonumber \\
\end{equation}
\begin{equation}
\frac{ck-2^n}{7}=+natural\hspace{0.1cm}  number-decimal \hspace{0.1cm} number  \rightarrow\frac{ck-2^n}{7}=decimal \hspace{0.1cm}  number \nonumber \\
\end{equation}
We reached a contradiction to the assumption. Thus, the theorem was verified. 
\end{proof}
\begin{example}
If the largest prime number was 13. \\
Step 1: $a=13$ \\ 
Step 2: $b= 11$ \\
Step 3: $c=11\times 9\times 7 \times 5 \times 3 \times 1$ \\
Step 4: If we know that 9 is not a prime number, we can write: \\ $ c=11 \times 7 \times 5 \times 3 \times 1=1155 $ \\
Step 5: $ R=ck-2^n=1155k-2^n \rightarrow 13<1155k-2^n \leq 168 $ 

\begin{equation}
k=1 \rightarrow \frac{13+2^n}{1155}<1 \leq \frac{168+2^n}{1155}    \rightarrow min \hspace{0.1cm} n=10 \nonumber \\
\end{equation}
\begin{equation}
\{   \frac{13+2^n}{1155}<k \leq \frac{168+2^n}{1155} \hspace{0.2cm}  (\frac{k}{2} \notin N ) \}  \nonumber \\
\end{equation}
\begin{equation} 
n=10 \rightarrow \frac{13+2^{10}}{1155}<k \leq  \frac{168+2^{10}}{1155}  (\frac{k}{2} \notin N) \rightarrow 0.89<k \leq 1.03 \rightarrow k=1  \nonumber \\
\end{equation}
\begin{equation}
n=11 \rightarrow \frac{13+2^{11}}{1155}<k \leq \frac{168+2^{11}}{1155}   (\frac{k}{2} \notin N) \rightarrow 1.78<k \leq 1.91\rightarrow k=- \nonumber \\
\end{equation}
\begin{equation}
\ldots  \nonumber \\
\end{equation}
\begin{equation}
n=18 \rightarrow \frac{13+2^{18}}{1155}<k \leq \frac{168+2^{18}}{1155} (\frac{k}{2} \notin N)  \rightarrow 226.97<k \leq 227.1 \rightarrow k=227  \nonumber \\
\end{equation}
\begin{equation}
n=19 \rightarrow \frac{13+2^{19}}{1155}<k \leq \frac{168+2^{19}}{1155}   (\frac{k}{2} \notin N)   \nonumber \\
\end{equation}
\begin{equation}
 \rightarrow 453.94<k \leq 454.07  \rightarrow  k=454,(\frac{k}{2} \notin N) \rightarrow  k=- \nonumber \\
\end{equation}
\begin{center}
And continuing so $\ldots$
\end{center}
Step 6: $ (k=1,n=10) \rightarrow R=131 ,  (k=227,n=18) \rightarrow R=41 $
\end{example}

\end{document}